\documentclass{fic-l}

\newtheorem{theorem}{Theorem}[section]
\newtheorem{lemma}[theorem]{Lemma}
\newtheorem{corollary}[theorem]{Corollary}

\theoremstyle{definition}

\theoremstyle{remark}
\newtheorem{remark}[theorem]{Remark}

\numberwithin{equation}{section}

\usepackage{amscd}

\begin{document}
\setcounter{page}{235}
\title{An example of symplectic 4-manifolds with positive signature}

\author{Martin Niepel}
\address{Department of Mathematics and Statistics\\ McMaster University\\ 
Hamilton\\ Ontario L8S 4K1}

\subjclass[2000]{Primary 57R17, 57R57; Secondary 53D35, 57N13}
\date{September 29, 2004}

\begin{abstract}
We construct a new family $\{K_n\}$ of simply connected symplectic 4-manifolds with the
property $c_1^2(K_n)/\chi_h(K_N) \rightarrow 9$ (as $n \rightarrow
\infty$).
\end{abstract}

\maketitle

\section{Introduction}
The homeomorphism type of a compact, closed, simply connected symplectic 4-manifold is 
uniquely determined by its Euler characteristic $e$, signature $\sigma$,
and the intersection form $Q$. If $Q$ is indefinite, as it happens in
most of the cases, it is determined by $e$,$\sigma$ and its \textit{parity}. It is natural to ask for 
which combinations ($e$, $\sigma$, $spin$) there exists such symplectic 
manifold, and how many distinct representatives there are. The first 
question is called the {\it geography problem}, the second is referred 
to as the {\it botany}. Both have been previously studied in the context
of complex surfaces.

Symplectic 4-manifolds can be equipped with an almost complex structure, 
and this provides the possibility to use some notions coming from the complex 
geometry. In particular, we can define Chern classes and Chern
numbers. Chern classes, in general, depend on the choice of an almost
complex structure, but $c_2$ is the same as Euler 
characteristic $e$, and $c_1^2$, which in the complex case denotes the 
self-intersection of a canonical divisor, can be extended uniquely by 
Hirzebruch signature formula $c_1^2=3 \sigma + 2e$. Similarly, we can
define the holomorphic Euler characteristic $\chi_h$ by setting 
$\chi_h = \frac{1}{4} (\sigma + e)$. The existence of the almost complex 
structure guarantees this number to be an integer.

The convention in 4-manifold geography is to use pairs ($\chi_h$, $c_1^2$) 
instead of ($e$, $\sigma$), and we will follow it in this article.

Important fact about Chern numbers is the 
Bogomolov--Miyaoka--Yau inequality, asserting that $c_1^2 \le 9 \chi_h$ 
holds for complex surfaces of general type. It was conjectured by 
by Fintushel and Stern that this relation holds also
for symplectic 4-manifolds. Most recent attempts (\cite{S1}, \cite{S2}, \cite{P})
in exploring the geography problem for simply connected symplectic 4-manifolds 
of positive signature, or equivalently those with $c_1^2 \ge 8 \chi_h$, provided
families of examples close to the Bogomolov--Miyaoka--Yau line (B--M--Y line).

These constructions followed similar pattern -- started with a complex
surface with nontrivial fundamental group lying on, or close to the B--M--Y 
line, and performed an operation which preserved symplectic structure and 
reduced the topological complexity of the original surface. The operation in 
use was either symplectic connected normal sum used several times, fibered
knot surgery or logarithmic transformation in a cusp neighborhood.

The botany question for complex surfaces has a partial answer -- the number 
of non-diffeomorphic complex surfaces of general type in a given homeomorphism type
is always finite, although there is no uniform bound on this
number. The other large family,  elliptic surfaces, is also well understood 
(see \cite{GS} for further references). However, for simply connected elliptic 
surfaces in a given homeomorphism type there are infinitely many
different smooth structures. These two results reflects the nature of 
restrictions given by the complex (geometric) structure.

In symplectic case, the picture is less clear due to the availability of
many gluing constructions, which are topological in nature. 
Fintushel and Stern in \cite{FS} provided another infinite family of symplectic 
4-manifolds with different Seiberg--Witten ($SW$) invariants, all of them 
homeomorphic to the standard $K3$ surface $E(2)$. The method, fibered knot surgery,
is different from the logarithmic transformation, and can be used
with similar success for any symplectic 4-manifold containing cusp fiber and
satisfying some additional properties (e.g. existence of a $c$-embedded
torus fiber). In particular, it can be used also for symplectic 
manifolds coming from (some) complex surfaces of general type. This indicates that 
the botany problem for symplectic 4-manifolds comprises infinitely many diffeomorphism 
types for a much larger class of homeomorphism types labeled by triples 
($\chi_h$, $c_1^2$, $spin$).

In the following we will construct a family $K_n$ of simply connected symplectic
4-manifolds close to the B--M--Y line. Moreover, for each $K_n$ there exist 
infinitely many mutually non-diffeomorphic symplectic 4-manifolds, all homeomorphic 
to $K_n$.

\begin{theorem} \label{thm1}
There are simply connected symplectic 4-manifolds $K_n$ with 
Chern numbers $c_2(K_n)= n^7 + 12 n^5 - 12n^4 + 6 n^3 + 22 $, 
$c_1^2(K_n)= 3n^7 + 20n^5 - 24 n^4 + 6n^3 +2$, holomorphic Euler characteristic  
$\chi_h (K_n)= \frac{1}{3}(n^7 + 8 n^5) - 3n^4 + n^3 +2$, and signature
$\sigma(K_n) = \frac{1}{3}(n^7 - 4n^5) - 2n^3 - 14$.
\end{theorem}

\textbf{Acknowledgement:} The author would like to thank Zolt\'an
Szab\'o for many helpful discussions and suggestions.

\section{Construction and homological properties of $X_n$}

In \cite{H}, Hirzebruch gave an example of a sequence $X_n$ of 
minimal complex surfaces of general type with
the following Chern numbers

$$ c_2(X_n) = n^7 $$
$$ c_1^2(X_n) = 3n^7 - 4n^5 .$$

\noindent
We briefly recall his construction, and also discuss some
topological properties of smooth models of these algebraic 
surfaces.

Let $\zeta = e^{2 \pi i /6}$ and $T$ be the elliptic curve

$$ T = \mathbb C /{ \{\mathbb Z \cdot 1 + \mathbb Z \zeta \} } . $$

\noindent
Let $T \times T$ be a complex surface, and denote its points by $(z,w)$. 
There is a well defined complex multiplication  
by algebraic integers $\alpha \in \mathbb Q [ \zeta ] $ on $T$. 
The elliptic curve in $T \times  T$ given by $w=\alpha z$ is then 
isomorphic to $T$ for every algebraic integer $\alpha$. 
Consider the following curves on $T \times T$ 

$$\begin{array}{ccc} T_0 :  w=0, & \quad &T_\infty :  z=0, \\
              T_1 : \; w=z, &  & T_\zeta :  w=\zeta z .
\end{array}$$

All four curves pass through the origin and do not intersect 
each other anywhere else. It is, in fact, the consequence of the 
special choice of $\zeta$ as a generator of the elliptic curve $T$.
Then $1-\zeta = - \zeta ^2$ is also a unit in $\mathbb Q[\zeta]$ and 
the curves $T_1$, $T_\zeta$ do not intersect outside of the origin.
Other cases are trivial.

Let $U_n$ be the lattice in $T \times T$ consisting of $n^4$ 
points

$$ U_n = \{ (z,w)|(nz,nw)=(0,0)\} .$$

For each point of $U_n$ there are four curves passing through
it, parallel to the curves $T_0$, $T_\infty$, $T_1$ and $T_\zeta$. 
Denote the union of $n^2$ curves parallel to $T_i$ as $D_i$ 
for $i \in \{0, \infty, 1, \zeta \}$.
Altogether we have $4n^2$ elliptic curves forming four parallel families. 
Except for the points in $U_n$ there are no other intersection 
points. After blowing up $n^4$ points of $U_n$ we get a smooth 
4-manifold $Y_n = T^4 \# n^4 \overline{\mathbb{CP}^2}$ with Euler characteristic $n^4$. 
There are $n^4$ exceptional 
curves $L_j$ ($j \in U_n$) resulting from blow-ups. Denote by  
$\tilde D_i$ the  proper transforms of $D_i$ after the blow-up; 
these curves are smooth and do not intersect 
each other.

Since each of the homology classes $[\tilde D_i]$ for $i \in \{0,\infty, 1, \zeta\}$ is 
divisible by $n^2$, we can construct a complex algebraic surface $X_n$
which is a branched $n^3$-fold cover of $Y_n$ branched over $\tilde D_i$ 
for each $i \in \{0,\infty, 1, \zeta\}$, where each $\tilde D_i$ 
will have the branching index $n$. Denote
this covering map as $\pi : X_n \rightarrow Y_n$.

Then a calculation of Euler characteristic gives

$$c_2(X_n)=n^3 \cdot e(Y_n \setminus \cup \tilde D_i) + n^2 \cdot 
e(\cup \tilde D_i)= n^3 \cdot n^4 = n^7 .$$

\noindent
For the Chern number $c_1^2$ we have

$$c_1^2(X_n)= 3n^7 - 4n^5 .$$

\noindent
For the actual computation, and for the further details 
we refer the reader to \cite{H}.

The map $\pi : X_n \rightarrow Y_n$ can be extended via the 
natural projections to the factors of the torus $T \times T$. 
Then we have the following diagram:

$$\begin{CD} @. X_n @. \\ @. @VV\pi V @. \\ T_0 @<p_1<< Y_n @>p_2>> T_\infty \end{CD}$$

The map $p_1$ respects the natural fibration of 
$T \times T$, singular fibers occur over the $n^2$ points in the set 
$U_{n,0}= U_n \cap T_0$, and are formed by the curve $\tilde D_\infty $ 
and exceptional curves $L_j$. The situation for the map $p_2$ 
is analogous.

Composed maps $p_1 \circ \pi$, $p_2 \circ \pi$ give two different 
(singular) fibrations of space $X_n$ over a 2-dimensional torus. In fact,  
these maps are holomorphic outside the blown-up points and corresponding 
exceptional curves, which leads us to the observation that a generic 
fiber of $p_1 \circ \pi$ is a smooth holomorphic curve. 

For a regular point
$z$ in $T_0$, i.e. point not belonging to $U_{n,0}$, the inverse image
$p_1^{-1}(z)$ is an elliptic curve parallel to $T_\infty$. 
This fiber is intersected $3n^2$ times by the curves $\tilde D_0, \tilde D_\zeta$ 
and $\tilde D_1$.
Note, that these points are branching points of the map $\pi$ restricted 
to the map of two--dimensional fibers in $X_n$ and $Y_n$ over the point $z$.
By the following computation of Euler characteristic of a regular fiber
$F_{\text{reg}} = (p_1 \circ \pi)^{-1}(z)$ we get
$$e(F_{\text{reg}}) = n^3\! \cdot e(T \setminus (3n^2 \: \text{pts.}))
+ n^2\! \cdot e(3n^2 \: \text{pts.}) =
n^3\! \cdot (-3n^2) + n^2\! \cdot 3n^2 = -3n^5 + 3n^4  .$$

\noindent
Hence the genus of a regular fiber is 
$g(F_{\text{reg}}) = 1 + \frac{3}{2} n^5 - \frac{3}{2} n^4$.

For the singular value $y$ of the composed map 
$p_1 \circ \pi$ and the respective singular fiber 
$F_{\text{sing}}=(p_1 \circ \pi)^{-1}(y)$ the situation is a bit 
more complicated. 

We know that $y$ belongs to $U_{n,0}$, and since the construction is 
completely symmetric with respect to the points in $U_{n,0}$, 
it is enough to describe the singular fiber over one point. We can choose 
the point $y=(0,0)$. 
Then the inverse image $p_1^{-1}(y)$ consists of $n^2$ exceptional 
spheres $L_j$ corresponding to the blown up points over 
$y$, connected with and intersected by the proper transform 
$\tilde T_\infty$ of the torus $T_\infty$.
Moreover, each of the exceptional spheres is intersected 
at exactly 3 different points by $\tilde D_0$, $\tilde D_\zeta$ and 
$\tilde D_1$, respectively. This gives us the complete information
about the map 
$\pi|_y : F_{\text{sing}} \rightarrow p_1^{-1}(y)$.

The singular fiber $F_{\text{sing}}$ consists of $n^2$ tori -- copies of 
$\tilde T_\infty$, connected with $n^2$ copies of branched covers 
$\tilde L_j$ of exceptional spheres $L_j$ ($j \in U_n \cap F_{\text{sing}}$). 
The Euler characteristic of $\tilde L_j$ is

$$e(\tilde L_j)= n^3 \cdot e(L_j \setminus 4\,\text{pts.}) + n^2 \cdot e( 4 \, \text{pts.}) =
-2 n^3 + 4 n^2 ,$$   

\noindent
and the Euler characteristic of the whole fiber $F_{\text{sing}}$ is

$$e(F_{\text{sing}})= n^2 \cdot e(\tilde L_j) + n^2 \cdot e(\tilde T_\infty) - n^4 \cdot e( \text{pt.})=
-2 n^5 + 3 n^4 .$$

To summarize the above mentioned facts we state the lemma.
 
\begin{lemma} \label{tor1}
There exists a complex surface $X_n$ of general type with Chern numbers 
$c_2(X_n)= n^7$,  $c_1^2(X_n)= 3 n^7 - 4 n^5$, holomorphic Euler characteristic 
$\chi_h(X_n)= \frac{1}{3} (n^7 - n^5)$, signature $\sigma(X_n) = \frac{1}{3}(n^7 - 4n^5)$,
and is equipped with a map
$p_1 \circ \pi : X_n \rightarrow T^2$. The genus of a regular fiber
$F_{\text{\rm{reg}}}$ is 
$g(F_{\text{\rm{reg}}})=  1 + \frac{3}{2} n^5 - \frac{3}{2} n^4$.
\qed \end{lemma}

\section{Topology of $X_n$}

In this section we describe the topology 
of $X_n$ via the study of its fundamental group. Since our
ultimate goal is to construct simply connected 4-manifold 
obtained from $X_n$ by a certain gluing, it would be 
convenient to find an embedded surface with small genus 
which generates as much of the fundamental group 
of $X_n$ as possible. We use the fact 
that the maps $p_1 \circ \pi$, $p_2 \circ \pi$ map the space 
$X_n$ to tori $T_0$ and $T_\infty$. 

Denote by $F_1 = (p_1 \circ \pi)^{-1} (z)$ a regular fiber over some 
regular point $z \in T_0$, and $F_2 = (p_2 \circ \pi)^{-1} (w)$ 
a regular fiber over some regular point $w \in T_\infty$. Recall 
that the resulting surfaces $F_1$, $F_2$ are smooth complex
curves. In the previous section we computed the genera of these curves, 
they are the same, and equal $g(F_1)=g(F_2)=  1 + \frac{3}{2} n^5 - \frac{3}{2} n^4$. 

The number of intersection points between $F_1$ and $F_2$ can be 
deduced from the following observation. Surface $F_2$ projects down
by the map $\pi$ to an embedded torus in $Y_n$. The map $\pi$ restricted 
to $F_2$ is then an $n^3$-fold branched cover. When we project resulting 
torus further down by the map $p_1$ into $T_0$ we get, actually, 
one-to-one map. Therefore the restricted map 
$p_1 \circ \pi|_{F_2} : F_2 \rightarrow T_0$ is an $n^3$-fold branched
cover as well. As a consequence we get that the number of intersection points 
between $F_1$ and $F_2$ is exactly $n^3$.
Note that all intersections are transversal and positive.

The following lemma shows that the union $F_1 \cup F_2$ contains 
all the topological information we need.

\begin{lemma}\label{fund1} The natural inclusion of fundamental groups
$i: \pi_1(F_1 \cup F_2) \rightarrow \pi_1(X_n)$ is surjective.
\end{lemma}

\begin{proof}
We will show that any loop $\gamma \in \pi_1(X_n)$ can be deformed 
into a loop lying inside $F_1 \cup F_2$. We can assume that the base 
point $p$ of the fundamental group $\pi_1 (X_n)$ is one of 
intersection points of $F_1$ and $F_2$.

The map $p_1 \circ \pi$ is singular over points in the set $U_{n,0}$.
Since the singular locus $(p_1 \circ \pi)^{-1} (U_{n,0})$ in $X_n$ 
has dimension 2, by using standard transversality argument we can 
assume that the representative of $\gamma$ avoids singularities.

Consider the projection $l=(p_1 \circ \pi)(\gamma)$. This gives a loop inside 
$T_0$ avoiding singular points in $U_{n,0}$. 
The inverse image of $l$, the space $(p_1 \circ \pi)^{-1}(l)$ is therefore
a surface bundle over $l$, and the fiber of this bundle over 
the point $p$ is exactly $F_1$. 
Using the fact that $F_2$ is a branched cover of $T_0$, we can deduce that 
the intersection of the bundle and the surface $F_2$ is an $n^3$-fold cover of
the loop $l$, in particular $F_2$ intersects each fiber in $n^3$ points.

Since the loop $\gamma$ lies inside the bundle whose fiber is connected, 
we can homotopically deform $\gamma$ into $\gamma \prime$ which will lie 
inside $F_2$ except for the part lying above the point $p$ and in the fiber $F_1$.
This proves the lemma.
\end{proof}

\begin{lemma} \label{tor3}
There exists an embedded symplectic surface $F$ in $X_n$ such that the 
inclusion map
$i : \pi_1 (F) \rightarrow \pi_1(X_n)$ is surjective. The surface 
$F$ has genus $g(F)= 3 n^5 - 3 n^4 + n^3 + 1$ and self-intersection $2n^3$.
\end{lemma}

\begin{proof}
We obtain $F$ from $F_1 \cup F_2$ by standard resolution of
singular intersection points in $F_1 \cap F_2$. It is easy to see
that $e(F)= 2 \cdot e(F_{\text{reg}}) - 2n^3 = -6 n^5 + 6 n^4 - 2n^3$.
Since $F_1$, $F_2$ are regular fibers, each of them have zero 
self-intersection number. The resolution of each positive 
intersection point increases the self-intersection by 2, hence 
$ [F]^2 = ([F_1] + [F_2])^2 = 2[F_1][F_2] = 2 n^3$.

Finally, the resolution of a singular point is only a local 
operation, but it might, and in fact does, change the fundamental
group of respective surfaces. Nevertheless, the map 
$\pi_1(F) \rightarrow \pi_1(F_1 \cup F_2)$ is surjective, its
kernel being generated by loops in $F$ which can be contracted
to singular points in $F_1 \cup F_2$. This proves the lemma.
\end{proof}

\section{Construction of $K_n$}

Surfaces $F_1$ and $F_2$ from the last section were non-degenerate
holomorphic curves. After the resolution of their intersection points 
we got a smooth embedded surface $F$, which is embedded symplectically.

We will perform connected symplectic sum operation (as defined by
Gompf in \cite{G})
to the (symplectic) manifold $X_n$ along embedded symplectic surface $F$
to obtain simply connected symplectic 4-manifold $K_n$. 
The ingredient we still need to provide is the second summand -- a simply 
connected symplectic manifold $N_n$, preferably with small Chern numbers $c_1^2$ 
and $c_2$, and with embedded symplectic surface of same genus and opposite
self-intersection as $F$. 

We begin with a $K3$ surface $E(2)$, 
equipped with the standard elliptic fibration $\pi : E(2) \rightarrow
\mathbb{CP}^1$. 
In order to obtain a non-degenerate holomorphic 
curve (or symplectic surface) of relatively large negative 
self-intersection, we need to select a number of points on a curve, 
and blow them up. Resulting proper transform of the original curve 
will have the required negative self-intersection, decreasing by 1 after each 
blow-up. The adjunction formula for complex and symplectic manifolds 
implies that this is basically the only effective procedure how to obtain 
such surfaces.

Picking and blowing-up $2n^3 -2$ points on a selected section of $E(2)$ gives  
a smoothly embedded rational complex curve $C$ with self-intersection $-2n^3$ inside 
complex surface $E(2) \# (2n^3 -2)\overline{\mathbb{CP}^2}$. 
Note that the elliptic fibration remained untouched outside blown-up points.
Now we can perform fibered knot surgery (see \cite{FS}), 
in a small tubular neighborhood of a regular elliptic fiber $T$ in a neighborhood 
of one of the cusp fibers in $E(2)$. Let $K \hookrightarrow S^3$ be a torus knot $K_{(p,q)}$ with genus $g=g(F)=\frac{1}{2}(p-1)(q-1)$, 
and $N_n=(E(2)\# (2n^3-2)\overline{\mathbb{CP}^2})_K$ be a result of fibered knot surgery on $E(2)\# (2n^3-2)\overline{\mathbb{CP}^2}$

\begin{lemma} \label{tor4}
The 4-manifold $N_n$ is symplectic and simply connected, and it contains symplectically 
embedded surface $F' = C \# \Sigma_g$ of genus $g=g(F)$ and self-intersection 
$[F']^2 = -2n^3 = - [F]^2$. Chern numbers of $N_n$ are $c_2(N_n) = 2n^3 + 22$,
$c_1^2(N_n)= -2n^3 + 2$, holomorphic Euler characteristic 
$\chi_h(N_n)= 2$, and the signature $\sigma(N_n)= -2n^3 - 14$. Moreover, 
$\pi_1(N_n \setminus F') = 1$.
\end{lemma}

\begin{proof}
The fact that $\pi_1(N_n \setminus F') = 1$ follows from the
triviality of $\pi_1(E(2) \setminus (C \cup T))$, since 
this space fibers over $S^2 \setminus \{\text{pt.}\}$ with a (singular) 
simply connected fiber $S^2 \setminus \{\text{pt.}\}$ coming from one of cusps. 
The rest of the statement follows directly from the construction of $N_n$ and straightforward computations.
\end{proof}

Now we can define $K_n$ as $K_n = X_n \#_F N_n$, and prove Theorem \ref{thm1}.

\begin{proof}\textbf{of Theorem \ref{thm1}}
We have done almost all the work in Lemmas \ref{tor1}, \ref{fund1},
\ref{tor3} and \ref{tor4}. 
It remains to show that the manifold $K_n$ is simply connected. This can be obtained by 
combining facts that $\pi_1(F) \rightarrow \pi(X_n)$ is surjective, that $\pi_1(N_n \setminus F')=1$,
and by using Seifert-Van Kampen theorem for the union $K_n= X_n \#_F N_n$. Computations of 
characteristics are routine.
\end{proof}

\begin{corollary}\label{infA} For each 4-manifold $K_n$ there 
exist infinitely many symplectic smooth manifolds of the same 
homeomorphism type which are not diffeomorphic to it. There 
exist infinitely many smooth manifolds of the same homeomorphism
type not admitting any symplectic structure, as well.
\end{corollary} 

\begin{proof}
The $K3$ surface $E(2)$
contains the nucleus $N(2)$, and at least one symplectically 
$c$-embedded torus in the complement $E(2) \setminus
N(2)$, because this complement is a \textit{Milnor fiber} $M(2,3,11)$ (see
\cite{GS}, \cite{S0}).
In our construction we used only the inside of the nucleus -- 
we used a regular fiber in the cusp neighborhood for the fibered
knot surgery $N_n=(E(2)\# (2n^3-2)\overline{\mathbb{CP}^2})_K$, 
and a section for the symplectic connected sum  $K_n = X_n \#_F N_n$. 
It means that the manifold $K_n$ could be produced in such a way 
that the $c$-embedded torus in the complement $E(2) \setminus N(2)$ 
stayed intact during the surgeries. This gives us a possibility to 
perform another fibered knot surgery in the neighborhood of the 
remaining cusp.

By reasoning similar to that in \cite{FS} we see that the manifolds 
after the knot surgery will have different $SW$-invariants than those 
of $K_n$, provided the Alexander polynomial of the knot is
non-trivial. The infinite number of fibered knots with different 
Alexander polynomials provide the infinite number of symplectic
manifolds homeomorphic but not diffeomorphic to $K_n$. It can be shown 
(see \cite{FS}, Cor. 1.3), that taking the fibered knot surgery along 
any knot, which is not fibered, and does not have a {\it monic}
Alexander polynomial will produce a manifold not admitting a
symplectic  structure. Again, there is an infinite number of such 
knots with different non-monic Alexander polynomials.
\end{proof}

\begin{remark} Similar results were obtained previously by J. Park, 
\cite{P}, and Stipsicz, \cite{S1}.
Both constructed infinite families of simply connected 4-manifolds 
with the ratio $c_1^2 / \chi_h$ very close to 9. However, their 
examples did not get above the line $c_1^2 = \alpha \cdot \chi_h$, 
with $\alpha \approx 8.76$, and $\alpha = 116/13 \approx  8.923$, respectively.

The 4-manifolds constructed by J. Park are {\it spin}. In this case 
the constructions are more rare, since there are fewer building blocks 
-- spin complex surfaces near B--M--Y line.

In \cite{S2}, Stipsicz improved his result and constructed a family of 
simply connected symplectic 4-manifolds with $(\chi_h, c_1^2)$
asymptotically approaching the line $c_1^2 = 9 \chi_h$. In comparison 
with our result, this family converges more rapidly, meaning that for 
comparable values of $\chi_h$ manifolds constructed by Stipsicz have 
higher value of $c_1^2$, and are closer to the B--M--Y line. Reason 
for this is, basically, the starting manifold. Our construction starts with 
$X_n$ below the line, while the other construction uses complex
surfaces lying on the B--M--Y line.

In \cite{park2}, J. Park used Stipsicz's method and populated a large
part of the region corresponding to the positive signature, however
the Euler characteristics of these 4-manifolds have to be quite large.

Some desirable results going further in this direction would be the 
construction of infinite family of 4-manifolds in constant distance 
from the B--M--Y line (compare with the results of Chen, \cite{Ch}, 
for complex surfaces of general type), or the construction of an example 
above it, which would contradict the conjecture by Fintushel and
Stern.  At the present time, the available tools and techniques are 
insufficient for such constructions.

Finally, we have to add that even the smallest manifolds constructed 
here are still rather large. The two smallest $K_n$'s with positive 
signature have the following invariants
$$ K_3: \qquad  \chi_h = 1 163, \quad c_1^2 = 9 641, \quad c_2 = 4 315, \quad \sigma = 227,$$
$$ K_4: \qquad  \chi_h = 7 490, \quad c_1^2 = 63 874, \quad c_2 = 26 006, \quad \sigma = 3 954.$$
\end{remark}


\end{document}